\documentclass[12pt]{amsart}

\usepackage{graphicx}
\usepackage{amsmath}
\usepackage{amscd}
\usepackage{amsfonts}
\usepackage{amssymb}
\usepackage{color}

\usepackage[unicode=true,
 bookmarks=true,bookmarksnumbered=false,bookmarksopen=false,
 breaklinks=false,pdfborder={0 0 1},backref=false,colorlinks=false]
 {hyperref}
\hypersetup{pdftitle={Essential normality, essential norms and hyperrigidity},
 pdfauthor={Matthew Kennedy and Orr Shalit}}

\numberwithin{equation}{section}

\newtheorem{theorem}{Theorem}[section]

\theoremstyle{definition}

\newcommand{\be}{\begin{equation}}
\newcommand{\ee}{\end{equation}}
\newcommand{\bes}{\begin{equation*}}
\newcommand{\ees}{\end{equation*}}

\newcommand{\bC}{\mathbb{C}}

\newcommand{\cK}{\mathcal{K}}

\newcommand{\cS}{\mathcal{S}}
\newcommand{\cT}{\mathcal{T}}

\newcommand{\cZ}{\mathcal{Z}}

\newcommand{\mb}[1]{\mathbb{#1}}

\newcommand{\pV}{\partial V}

\begin{document}

\title[Corrigendum]{Corrigendum to: Essential normality, essential norms and hyperrigidity}

\author{Matthew Kennedy} 
\address{Department of Pure Mathematics \\
University of Waterloo \\
200 University Avenue West \\
Waterloo, Ontario N2L 3G1 \\
Canada}
\email{matt.kennedy@uwaterloo.ca}

\author[Orr Shalit]{Orr Moshe Shalit} 
\address{Faculty of Mathematics\\
Technion - Israel Institute of Technology\\
Haifa\; 12000\\
Israel}
\email{oshalit@tx.technion.ac.il}

\thanks{First author supported by grant no. 418585 from the Natural Sciences and Engineering Research Council of Canada. 
The second author is partially supported by ISF Grant 474/12 and  EU FP7/2007-2013 Grant 321749. }

\begin{abstract} 
In our paper ``Essential normality, essential norms and hyperrigidity'' we claimed that the restriction of the identity representation of a certain operator system (constructed from a polynomial ideal) has the unique extension property, however the justification we gave was insufficient. 
In this note we provide the required justification under some additional assumptions. Fortunately, homogeneous ideals that are ``sufficiently non-trivial'' are covered by these assumptions. This affects the section of our paper relating essential normality and hyperrigidity. We show here that Proposition 4.11 and Theorem 4.12 hold under the additional assumptions. We do not know if they hold in the generality considered in our paper.
\end{abstract}


\maketitle

\section{Statement of the problem and consequences} 
We use the notation and terminology from \cite{KenSha15}.  In Proposition 4.11 and Theorem 4.12 of \cite{KenSha15}, we claimed that for a polynomial ideal $I \triangleleft \bC[z_1,\ldots,z_d]$ satisfying certain standing assumptions, the restriction of the identity representation of the C*-algebra $\cT_{I}$ to the operator system $\cS$ has the uniqe extension property. As justification we cited \cite[Proposition 6.4.6]{ChenGuo}, which treats submodules but not quotient modules. We believed that the same proof works, and omitted further detail.

In fact, the same proof does work in the case of homogeneous ideals that are ``sufficiently non-trivial'' (which is the case of primary interest for the conjecture of Arveson that our paper treats). However, it does not work in the generality that we required. If the ideal is not ``sufficiently non-trivial,'' then not only does the proof break down, the result is false: the identity representation is {\em not} a boundary representation. We are grateful to Michael Hartz and Raphael Clouatre who pointed out these issues. 

The standing assumptions we imposed on $I$ were that $d \geq 2$, that
\be\label{eq:stand_assum}
\overline{V}(I) := \overline{\cZ(I) \cap \mb{B}_d} = \cZ(I) \cap \overline{\mb{B}_d},
\ee
and that
\be\label{eq:stand_assum2}
\pV := \cZ(I) \cap \partial  \mb{B}_d \neq \emptyset , 
\ee
where $\cZ(I)$ is the variety determined by $I$. 
These assumptions ensure that $I$ does not have finite codimension in $\mb{C}[z] = \mb{C}[z_1, \ldots, z_d]$. 
The only reason for the $d \geq 2$ assumption is to ensure that we do not consider the zero ideal in $\mb{C}[z_1]$. For, in that case, the corresponding quotient module is the Hardy space $H^2(\mb{D})$, which is well understood and has radically different properties. 
For example, the operator system generated by the unilateral shift on $H^2(\mb{D})$ has a non-trivial Shilov boundary in the Toeplitz algebra, meaning in particular that the identity representation does not have the unique extension property. 

In general, the assumptions \eqref{eq:stand_assum} and \eqref{eq:stand_assum2} are not enough to ensure that the identity representation has the unique extension property. 
We need to assume in addition that
\be\label{eq:stand_assum3}
\textrm{$0 \in \cZ(I)$ and $0$ is a not isolated in $\cZ(I)$},
\ee
and that
\be\label{eq:stand_assum4}
\textrm{the ideal $I$ contains no linear polynomials}. 
\ee
Assumption \eqref{eq:stand_assum4} is a harmless non-trivialtiy assumption. 
Indeed, if $I$ contains a linear polynomial then we may assume that it is $z_d$, and in this case we are essentially working in a lower dimensional situation. 
On the other hand, it is necessary to assume something along these lines, since if $d=2$ and $I$ is the ideal generated by $z_2$, then the quotient module $H^2_2 / I$ is precisely the Hardy module $H^2(\mb{D})$, and as noted above, the unique extension property fails.

We will say that the ideal $I$ is {\em sufficiently non-trivial} when $I \ne 0$, $I$ is not of finite codimension and assumption \eqref{eq:stand_assum4} holds. 

In constrast, assumption \eqref{eq:stand_assum3} does not seem so harmless, and we do not know whether it is necessary. 
However, all of the assumptions \eqref{eq:stand_assum}-\eqref{eq:stand_assum4} hold for any homogeneous (or even quasi-homogeneous) ideal that is sufficiently non-trivial, and as mentioned, this is the case of primary interest for Arveson's conjecture.

We now state the main result of this note which fills the gap in \cite{KenSha15}. 

\begin{theorem}\label{thm:main}
Let $I \triangleleft \bC[z_1,\ldots,z_d]$ be a polynomial ideal satisfying the assumptions \eqref{eq:stand_assum}, \eqref{eq:stand_assum2}, \eqref{eq:stand_assum3} and \eqref{eq:stand_assum4}. Then the identity representation of $\cT_I$ is a boundary representation for $\cS$.  In particular, the restriction of the identity representation of $\cT_I$ to $\cS$ has the unique extension property. 
\end{theorem}

We note that Proposition 4.11 and Theorem 4.12 from our paper \cite{KenSha15} are only known to hold under the assumptions in Theorem \ref{thm:main}. In particular, these results hold for all sufficiently non-trivial homogeneous ideals. However, we do not know what happens in general. The other results in our paper are not affected.

\section{Proof of Theorem \ref{thm:main}}
\label{sec:proofmain}

The algebra $\cT_I$ is irreducible and contains the compacts, hence by Arveson's boundary theorem \cite[Theorem 2.1.1]{Arv72} it suffices to show that the quotient map $q : \cT_I \to \cT_I / \cK$ is not completely isometric on $\cS$. 

Suppose that $I$ is generated as an ideal by polynomials $\{f_1, \ldots, f_k\}$. 

\noindent{\bf Claim:} There is a choice of basis such that not all of the variables $z_1, \ldots, z_d$ appear in the linear terms of $f_1, \ldots, f_k$.

\noindent{\bf Proof of claim:} If the degree one term of each $f_i$ is zero, then we are done. Hence we can assume that at least one of the $f_i$'s, say $f_1$, has a nonzero linear term. Notice that by assumption (\ref{eq:stand_assum3}), the constant term of $f_1$ is zero. Hence without loss of generality, we can suppose that
\[
f_1 = z_1 + \textrm{higher order terms} .
\]
If needed, we can subtract multiplies of $f_1$ from $f_2,\ldots,f_k$, and hence we can assume that $z_1$ does not appear in the linear term of each $f_2, \ldots, f_k$. 
If the linear term of each $f_2, \ldots, f_k$ is zero, or if $k = 1$ while $d > k = 1$, then we are done. Otherwise we can repeat this process a finite number of times until we find ourselves in the situation where necessarily $k \geq d$ and 
\[
f_i = z_i + \textrm{higher order terms} \,\, , \,\, i=1, \ldots, d. 
\]
Then the map $F: \mb{C}^d \to \mb{C}^k$ given by $F(z) = (f_1(z), \ldots, f_k(z))$ for $z \in \mb{C}^d$ satisfies $F(0) = 0$. Moreover, the derivative $dF(0)$ has full rank, and it follows that $0$ is an isolated point in $\cZ(I)$, contradicting assumption (\ref{eq:stand_assum3}). This proves the claim. 

Now assuming that $z_1$ does not appear in the linear terms of any of the generators of $I$, we conclude that $S_1^* S_1 1 = 1$, because $z_1 \in I^\perp$.
Furthermore, we also have $\|S_2^* S_2 1\| > 0$, because $z_2 \notin I$.
We conclude that the norm of the column operator $\left[ \begin{smallmatrix}  S_1\\  S_2\end{smallmatrix}\right]$ is strictly bigger than $1$. 
On the other hand, the norm of the row operator $\left[ \begin{smallmatrix}  Z_1\\  Z_2\end{smallmatrix}\right]$, which is the norm of $\left[ \begin{smallmatrix}  S_1\\  S_2\end{smallmatrix}\right]$ in the quotient by the compacts, is equal to $1$, thanks to Theorem 3.3 in \cite{KenSha15}. 
This completes the proof.

\section{Other spaces}
We conclude this note by pointing out that although the proof of Theorem \ref{thm:main} was worked out in the setting of quotients of $H^2_d$, a variant also works in the setting of the Besov-Sobolev spaces considered in \cite[Section 5]{KenSha15}. 
In fact, one aspect of the situation is slightly simpler here, since we only consider homogeneous ideals. In particular, \cite[Corollary 5.3 (4)]{KenSha15} still holds for all sufficiently non-trivial homogeneous ideals.

Let us spell out the details. Suppose $I$ is a sufficiently non-trivial homogenous ideal in $\mb{C}[z_1,\ldots,z_d]$. Fix $\sigma \in [1/2, d/2)$, and let $T$ be the compression of $M_z$ to $B^2_\sigma(\mb{B}_d) \ominus I$. 
Since $I$ is sufficiently non-trivial, it contains no linear terms. 
Therefore, invoking equation \cite[(5.2)]{KenSha15}, we find that $\sum_{i=1}^d T_i^* T_i 1 = \frac{d}{2\sigma}1$. Also, because $\sigma < d/2$ we find that the column operator  $\left[ \begin{smallmatrix}  T_1 & \cdots & T_d \end{smallmatrix}\right]^t$ has norm at least $\frac{d}{2\sigma} > 1$. 
On the other hand, the essential norm of $\left[ \begin{smallmatrix}  T_1 & \cdots & T_d \end{smallmatrix}\right]^t$ is $1$, by \cite[Corollary 5.3 (2)]{KenSha15}. 
Since $C^*(T)$ is an irreducible algebra containing the compacts, Arveson's boundary theorem implies that the identity representation of $\mathrm{C}^*(T_1,\ldots,T_d)$ is a boundary representation for $T$, so its restriction to the operator system generated by $T$ has the unique extension property.

\bibliographystyle{amsplain}

\end{document}